\documentclass[11pt]{amsart}

\usepackage{amsmath,amssymb,amsthm}
\usepackage{fancyhdr}
\usepackage{amsfonts,graphicx}

\theoremstyle{plain}
\newtheorem{Theo}{Theorem}[section]
\newtheorem{lem}[Theo]{Lemma}
\newtheorem{prop}[Theo]{Proposition}

\theoremstyle{plain}
\theoremstyle{definition}

\theoremstyle{remark}
\newtheorem{Rema}{\bf Remark}
\newtheorem*{rema*}{Remarks}

\newcommand{\ZZ}{\mathbb{Z}}  
  
\newcommand{\NN}{\mathbb{N}}

\newcommand{\RR}{\mathbb{R}} 
 
\newcommand{\Lip}{\textnormal{Lip}}

\newcommand{\diver}[1]{ \textnormal{div}\hspace{0.07cm} #1} 

\date{}

\begin{document}
\title{ On the global  well-posedness of the  Boussinesq system  with zero viscosity}
\author{Taoufik Hmidi \& Sahbi Keraani}

\begin{abstract}
In this paper we prove the global well-posedness of the two-dimensional Boussinesq system with zero viscosity  for rough initial data.

\end{abstract}
\maketitle
\section{Introduction}
This paper is a sequel to \cite{T}. We continue to study  the global existence for the 
two-dimensional  Boussinesq system,
  $$({B_{\nu,\kappa}})\left\lbrace
\begin{array}{l}
\partial_t v+v\cdot\nabla v- \nu\Delta
 v+\nabla \pi=\theta e_{2},\\
 \partial_{t}\theta+v\cdot\nabla\theta-\kappa\Delta \theta=0,\\
\textnormal{div}\hspace{0.07cm}v=0,\\
 {v}_{| t=0}=v ^0,\quad \theta_{|t=0}=\theta ^0.
\end{array}
\right.$$ 
 Here, $e_{2}$ denotes the vector $(0,1),$ $v=(v_{1},v_{2})$ is the velocity field, $\pi$ the scalar pressure and $\theta$ the temperature. The coefficients
  $\nu$ and $\kappa$ are assumed to be positive; $\nu$ is called the kinematic viscosity and $\kappa$ the molecular conductivity.\\
    In the case of strictly positive coefficients $\nu$ and $\kappa$ both velocity and temperature have sufficiently smoothing effects leading to the
     global well-posedness results  proven by numerous authors  in various function spaces 
  (see  \cite{Can, Guo,S} and the references therein).\\
  For $\nu>0$ and $\kappa=0$ the problem of global well-posedness is well understood. In \cite{Cha}, Chae proved global well-posedness for initial data $(v^0,\theta^0)$ lying in Sobolev spaces $H^s\times H^s,$ \mbox{with $s>2$} ( see also \cite{Hou}). This result has been recently improved in \cite{T0} by taking the data in $H^s\times H^s,$ with $s>0.$ However we give only a global existence result without uniqueness in the energy space $L^2\times L^2.$ In \cite{A} we prove a uniqueness result for data belonging to $ L^2\cap B_{\infty,1}^{-1}\times B_{2,1}^0.$ More recently Danchin and Paicu \cite{rd} have established a uniqueness result in the energy space.
    
  Our goal here is to study the global well-posedness of the system $(B_{0,\kappa}),$ with $\kappa>0.$   First of all, let us recall that the two-dimensional incompressible Euler system, corresponding to $\theta^0=0,$ 
is globally well-posed in the Sobolev \mbox{space $H^s,$} with $s>2.$ This is due to the  advection  of the vorticity by the flow: there is no accumulation of the vorticity  and thus there is no finite time singularities according to B-K-M criterion \cite{M}. In critical spaces like $B_{p,1}^{\frac{2}{p}+1}$ the situation is more complicate because we do not know if the B-K-M criterion works or not.  In \cite{v1}, 
  Vishik proved  that Euler system is globally well-posed in these
   critical Besov spaces.  He used for the proof a new logarithmic estimate taking advantage on the  particular structure of the vorticity equation in dimension two.    For the Boussinesq \mbox{system $(B_{0,\kappa})$},  Chae 
   has proved in
    \cite{Cha} the global well-posedness  for initial data $v^0,\theta^0$ lying in Soboloev space $H^s,$ \mbox{with $s>2,$}  His method is basically related to Sobolev logarithmic estimate in which the velocity and the temperature are needed to be Lipschitzian explaining the restriction $s>2.$   We intend here to improve this result for  rough initial data. Our results reads as follows (see the definition of Besov spaces
       given in next section ).
    \begin{Theo}\label{thmm0}
     Let   $v^0\in B_{p,1}^{1+\frac 2p}$ be a divergence-free vector field of $\RR^2$ and $\theta ^0\in L^r,$ with $2<r\leq p<\infty$ . Then there exists a unique global solution $(v,\theta)$ to the Boussinesq 
system $(B_{0,\kappa}), \kappa>0$ such that
$$
v\in \mathcal{C}(\RR_{+}; B_{p,1}^{1+\frac2p})\,\quad \hbox{and}\quad\theta\in L^\infty(\RR_{+};L^r)\cap \widetilde L^1_{\textnormal{loc}}(\RR_{+}; B_{r,\infty}^2).
$$
    \end{Theo}
    The situation in the case  $p=+\infty$ is more subtle since Leray's projector is not continuous on $L^\infty$ and we overcome this  by working in homogeneous Besov spaces leading to more technical difficulties. Before stating our result we introduce the following sub-space \mbox{of $L^\infty:$}
    $$
    u\in \mathcal{B}^\infty\Leftrightarrow \|u\|_{\mathcal{B}^\infty}:=\|u\|_{L^\infty}+\|\Delta_{-1}u\|_{\dot B_{\infty,1}^0}<\infty.
    $$
   We notice that $\mathcal{B}^\infty$ is a Banach space and independent of the choice of the unity dyadic partition.  For the definition of the frequency localization operator $\Delta_{-1}$ we can see next section.
    Our second main   result is the following:
 \begin{Theo}\label{thm1}
 Let   $v^0\in B_{\infty,1}^{1},$ with zero divergence  and $\theta ^0\in \mathcal{B}^\infty$  . Then there exists a unique global solution $(v,\theta)$ to the Boussinesq 
system $(B_{0,\kappa}), \kappa>0$ such that
$$
v\in \mathcal{C}(\RR_{+}; B_{\infty,1}^{1})\,\quad \hbox{and}\quad\theta\in L^\infty_{\textnormal{loc}}(\RR_{+}; \mathcal{B}^\infty)\cap \widetilde L^1_{\textnormal{loc}}(\RR_{+}; B_{\infty,\infty}^2).
$$
 
 \end{Theo}
  The key of the proof is to bound for every time Lipschitz norm of both velocity and  temperature. This will be done by using some logarithmic estimates and Theorem \ref{effect}.  The last one describes new smoothing effects for the transport-diffusion equation governed by a vector field which is not necessary Lipschitzian but only quasi-lipschitzian.  Its proof is done in the spirit of \cite{T0}. 
  
  \
   
 The rest of this paper is organized as follows. In section $ 2$, we recall some preliminary results on Besov spaces.  Section $ 3$ is devoted to the proof of smoothing effects. In section $4$ and $5$  we give respectively the proof of  {Theorem \ref{thmm0} and \ref{thm1}.}  We give in the appendix  a logarithmic estimate and a commutator lemma. 

\section{Notation and preliminaries} 
Throughout this paper  we shall denote by $C$ some real positive constants which may be different in each occurrence and by $C_{0}$ a real positive
 constant depending on the initial data.

Let us introduce the so-called   Littlewood-Paley decomposition and the corresponding  cut-off operators. 
There exists two radial positive  functions  $\chi\in \mathcal{D}(\RR^d)$ and  $\varphi\in\mathcal{D}(\RR^d\backslash{\{0\}})$ such that
\begin{itemize}
\item[\textnormal{i)}]
$\displaystyle{\chi(\xi)+\sum_{q\geq0}\varphi(2^{-q}\xi)=1}$;$\quad \displaystyle{\forall\,\,q\geq1,\, \textnormal{supp }\chi\cap \textnormal{supp }\varphi(2^{-q})=\varnothing}$\item[\textnormal{ii)}]
 $ \textnormal{supp }\varphi(2^{-p}\cdot)\cap
\textnormal{supp }\varphi(2^{-q}\cdot)=\varnothing,$ if  $|p-q|\geq 2$.
\end{itemize}

For every $v\in{\mathcal S}'(\RR^d)$ we set 
  $$
\Delta_{-1}v=\chi(\hbox{D})v~;\, \forall
 q\in\NN,\;\Delta_qv=\varphi(2^{-q}\hbox{D})v\quad\hbox{ and  }\;
 S_q=\sum_{-1\leq p\leq q-1}\Delta_{p}.$$
 The homogeneous operators are defined by
 $$
 \dot{\Delta}_{q}v=\varphi(2^{-q}\hbox{D})v,\quad \dot S_{q}v=\sum_{j\leq q-1}\dot\Delta_{j}v,\quad\forall q\in\ZZ.
 $$
From   \cite{b}  we split the product 
  $uv$ into three parts: $$
uv=T_u v+T_v u+R(u,v),
$$
with
$$T_u v=\sum_{q}S_{q-1}u\Delta_q v \quad\hbox{and}\quad R(u,v)=\sum_{|q'-q|\leq 1}\Delta_qu\Delta_{q'}v.
$$
  
Let us now define inhomogeneous and homogeneous Besov spaces. For
 $(p,r)\in[1,+\infty]^2$ and $s\in\RR$ we define  the inhomogeneous Besov \mbox{space $B_{p,r}^s$} as 
the set of tempered distributions $u$ such that
$$\|u\|_{B_{p,r}^s}:=\Big( 2^{qs}
\|\Delta_q u\|_{L^{p}\Big)_{\ell ^{r}}}<+\infty.$$
The homogeneous Besov space $\dot B_{p,r}^s$ is defined as the set of  $u\in\mathcal{S}'(\RR^d)$ up to polynomials such that
$$
\|u\|_{\dot B_{p,r}^s}:=\Big( 2^{qs}
\|\dot\Delta_q u\|_{L^{p}}\Big)_{\ell ^{r}(\ZZ)}<+\infty.
$$
{Let $T>0$} \mbox{and $\rho\geq1,$} we denote by $L^\rho_{T}B_{p,r}^s$ the space of distributions $u$ such that 
$$
\|u\|_{L^\rho_{T}B_{p,r}^s}:= \Big\|\Big( 2^{qs}
\|\Delta_q u\|_{L^p}\Big)_{\ell ^{r}}\Big\|_{L^\rho_{T}}<+\infty.$$
We say that 
$u$ belongs to the space
 $\widetilde L^\rho_{T}{B_{p,r}^s}$ if
 $$
 \|u\|_{ \widetilde L^\rho_{T}{B_{p,r}^s}}:= \Big( 2^{qs}
\|\Delta_q u\|_{L^\rho_{T}L^p}\Big)_{\ell ^{r}}<+\infty .$$
The relations between these spaces are detailed below are a direct consequence of 
the Minkowski inequality. 
%\begin{lem}\label{lemm4}
Let $ \epsilon>0,$ then 
$$
L^\rho_{T}B_{p,r}^s\hookrightarrow\widetilde L^\rho_{T}{B_{p,r}^s}\hookrightarrow{L^\rho_{T}}{B_{p,r}^{s-\epsilon}}
,\,\textnormal{if}\quad  r\geq \rho,$$
$$
{L^\rho_{T}}{B_{p,r}^{s+\epsilon}}\hookrightarrow\widetilde L^\rho_{T}{B_{p,r}^s}\hookrightarrow L^\rho_{T}B_{p,r}^s,\, \textnormal{if}\quad 
\rho\geq r.$$
%\end{lem}
 We will  make continuous use of Bernstein inequalities (see for \mbox{example \cite{Ch1}}).
\begin{lem}\label{lb}\;
 There exists a constant $C$ such that for $k\in\NN,$ $1\leq a\leq b$ and for  $u\in L^a(\RR^d)$, 
\begin{eqnarray*}
\sup_{|\alpha|=k}\|\partial ^{\alpha}S_{q}u\|_{L^b}&\leq& C^k\,2^{q(k+d(\frac{1}{a}-\frac{1}{b}))}\|S_{q}u\|_{L^a},\\
\ C^{-k}2^
{qk}\|\dot{\Delta}_{q}u\|_{L^a}&\leq&\sup_{|\alpha|=k}\|\partial ^{\alpha}\dot{\Delta}_{q}u\|_{L^a}\leq C^k2^{qk}\|\dot{\Delta}_{q}u\|_{L^a}.
\end{eqnarray*}

\end{lem}
The following result is due to Vishik \cite{v1}.
\begin{lem}\label{Vishik00}
Let $d\geq 2$, there exists a positive constant  $C$ such that for any smooth function $f$  and for any 
diffeomorphism  $\psi$ of $\RR^d$ preserving  Lebesgue measure, we have for all
$p\in[1,+\infty]$ and for all  { $j,q\in\ZZ,$}
$$\|\dot\Delta_j(\dot\Delta_q f\circ\psi)\|_{L^p}\leq
C2^{-|j-q|}\|\nabla\psi ^{\eta(j,q)}\|_{L^{\infty}}\|\dot\Delta_q
f\|_{L^p},$$
with
$$\eta(j,q)=\hbox{sign}(j-q).$$
\end{lem}
      Let us now recall the following  result proven in \cite{rd,T1}.
      \begin{prop}\label{propagation}
      Let $\nu\geq 0,\,(p,r)\in [1,\infty]^2,\,s\in]-1,1[, $ $v\in  L^1_{\textnormal{loc}}(\RR_{+};\textnormal{Lip}(\RR^d))$ with zero divergence and $f$ be a smooth fuction. Let $a$ be any smooth solution of the transport-diffusion equation
      $$
\partial_{t} a+v\cdot\nabla a-\nu\Delta a=f.      $$
  Then there is a constant $C:C(s,d)$ such that  for every $t\in\RR_{+}$
  $$
  \|a\|_{\widetilde L^\infty_{t}B_{p,r}^s}+\nu^{\frac1m}\|a-\Delta_{-1}a\|_{\widetilde L^m_{t}B_{p,r}^{s+\frac2m}}\le Ce^{CV(t)}\Big(\|a^0\|_{B_{p,r}^s}+\int_{0}^t\|f(\tau)\|_{B_{p,r}^s}d\tau\Big),
  $$
  where $V(t):=\displaystyle \int_{0}^t\|\nabla v(\tau)\|_{L^\infty}d\tau.$
  
 % Furthermore the above inequality holds ture for $s=-1, r=\infty$ and $p\in[1,\infty],$ despite we change $V(t)$ by $\|v\|_{L^1_{t}B_{\infty,1}^1}.$
  \end{prop}
   \section{Smoothing effects}
   This section is devoted to the proof of  a new regularization effect for a transport-diffusion equation with respect to a vector field which is not necessary Lipschitzian. This problem was studied by the first author \cite{T0} in the context of singular vortex patches for two dimensional Navier-Stokes equations. The estimate given below is more precise.
   \begin{Theo}\label{effect}
   Let  $v$ be  a smooth  divergence-free vector field of $\RR^d $ with vorticity $\omega:=\textnormal{curl } v.$ Let $a$ be a smooth solution of the transport-diffusion equation
   $$
   \partial_{t} a+v\cdot\nabla a-\Delta a=0;\quad a_{|t=0}=a^0.
   $$
   Then we have for $q\in\NN\cup\{ -1\}$ and $t\geq 0$
   \begin{eqnarray*}
   2^{2q}\int_{0}^t\|\Delta_{q}a(\tau)\|_{L^\infty}d\tau&\lesssim&\|a^0\|_{L^\infty}\Big(1+t+(q+2)\|\omega\|_{L^1_{t}L^\infty}
   +   \|\nabla\Delta_{-1}v\|_{L^1_{t}L^\infty}\Big).
 \end{eqnarray*}

   \end{Theo}
   \begin{Rema}
   In \cite{T1}, the first author proved in the case of Lipschitzian velocity the following estimate   \begin{equation}\label{visc0}
   2^{2q}\int_{0}^t\|\Delta_{q}a(\tau)\|_{L^\infty}d\tau\lesssim\|a^0\|_{L^\infty}\Big(1+ t+\int_{0}^t\|\nabla v(\tau)\|_{L^\infty}d\tau\Big).
   \end{equation}
   We emphasize that the above theorem is also true when we change $L^\infty$ by  $L^p,$ \mbox{with $p\in[1,\infty].$}
   \end{Rema}
   \begin{proof}
  The idea of the proof is the same as in \cite{T1}. We use Lagrangian formulation combined with intensive use of paradifferential calculus.

Let $q\in\NN^*,$ then the Fourier localized function  $a_{q}:=\Delta_{q} a$ satisfies

\begin{equation}
\label{Eq:1}
\partial_{t} a_{q}+S_{q-1}v\cdot\nabla a_{q}-\Delta a_{q}=(S_{q-1}-\hbox{Id})v\cdot\nabla a_{q}-[\Delta_{q},v\cdot\nabla] a:=g_{q}.
\end{equation}
Let $\psi_{q}$ denote the flow of the regularized velocity $S_{q-1}v$:
$$
\psi_{q}(t,x)=x+\int_{0}^tS_{q-1}v\big(\tau,\psi_{q}(\tau,x)\big)d\tau.
$$
We set
$$
\bar{a}_{q}(t,x)=a_{q}(t,\psi_{q}(t,x))\quad\hbox{and}\quad \bar{{g}}_{q}(t,x)={g}_{q}(t,\psi_{q}(t,x)).
$$
From Leibnitz formula we deduce the following identity
\begin{equation}
\label{14}
\nonumber \Delta \bar{a}_{q}(t,x)=   \sum_{i=1}^{d}\Big\langle H_{q}\cdot(\partial ^i \psi_q
)(t,x),(\partial ^i \psi_q )(t,x)\Big\rangle +  (\nabla a_{q})(t,\psi_{q}(t,x))\cdot\Delta\psi_{q}(t,x)   ,
\end{equation}
where $H_{q}(t,x):=
 (\nabla
^2{a}_q)(t,\psi_{q}(t,x))
$ is the Hessian matrix.\\
 Straightforward computations based on the definition of the flow and Gronwall's inequality yield
$$
\partial ^i\psi_{q} (t,x)=e_i+h_{q}^i(t,x),
$$ 
where $(e_{i})_{i=1}^d$ is the canonical basis of $\RR^d$ and the function $h_{q}^{i}$ is estimated as follows
\begin{equation}
\label{marha}
\|h_{q}^{i}(t)\|_{L^\infty}\lesssim V_{q}(t)e^{CV_{q}(t)},\quad\hbox{with}\quad V_{q}(t):=\int_{0}^t\|\nabla S_{q-1}v(\tau)\|_{L^\infty}d\tau.
\end{equation}
Applying Leibnitz formula and Bernstein inequality we find
\begin{equation}
\label{zada}
\|\Delta\psi_{q}(t)\|_{L^\infty}\lesssim 2^{q}V_{q}(t)e^{CV_{q}(t)}.
\end{equation}
The outcome is 
\begin{equation}
\label{Tr45}
\Delta \bar{a}_{q}(t,x)=(\Delta a_{q})(t,\psi_{q}(t,x))-\mathcal{R}_{q}(t,x),
\end{equation}
with
\begin{eqnarray}\label{Eq:3}\nonumber
\|\mathcal{R}_{q}(t)\|_{L^\infty}&\lesssim& \|\nabla a_{q}(t)\|_{L^\infty}\|\Delta \psi_{q}(t)\|_{L^\infty}\\
\nonumber&+&\|\nabla^2 a_{q}(t)\|_{L^\infty}\sup_{ i}\big(\|h^{i}_{q}(t)\|_{L^\infty}+\|h^{i}_{q}(t)\|_{L^\infty}^2\big)
\\ 
&\lesssim&
2^{2q}V_{q}(t)e^{CV_{q}(t)}\|a_{q}(t)\|_{L^\infty}.
\end{eqnarray}
In the last line we have used  Bernstein inequality.

From (\ref{Eq:1}) and (\ref{Tr45}) we see that $\bar{a}_{q}$ satisfies
\begin{equation}
\label{T199}
\nonumber (\partial_t-\Delta)\bar{a}_q(t,x)=\mathcal{R}_{q}(t,x)+\bar{g}_q(t,x).
\end{equation}
Now, we will again localize in frequency  this equation through the \mbox{operator $\Delta_{j}$. }So we write from  Duhamel formula,
\begin{eqnarray}
\label{ali}
\nonumber\Delta_j\bar{a}_q(t,x)&=&e ^{ t\Delta}\Delta_ja_q(0)+\int_{0}^{t}
e ^{ (t-\tau)\Delta}\Delta_j\mathcal{R}_q(\tau,x)d\tau\\
&+&\int_{0}^{t}
e ^{ (t-\tau)\Delta}\Delta_j\bar{g}_q(\tau,x)d\tau.
\end{eqnarray}
At this stage we need the following lemma (see for instance \cite{Ch2}).
\begin{lem}
For  $u\in L^\infty$ and  $j\in\NN,$
\begin{equation}
\label{mohima1}
\|e^{ t\Delta}\Delta_{j}u\|_{L^\infty}\leq C e^{-c t2^{2j}}\|\Delta_{j}u\|_{L^\infty},
\end{equation}
where the constants $C$ and $c$ depend only on the dimension $d$.
\end{lem}
Combined with (\ref{Eq:3})  this lemma yields, for every $j\in\NN$,
  \begin{equation}
\label{Eq:4}
\|e^{(t-\tau)\Delta}\Delta_{j}\mathcal{R}_{q}(\tau)\|_{L^\infty}
\lesssim 2^{2q}V_{q}(\tau)e^{CV_{q}(\tau)} e^{-c(t-\tau)2^{2j}}\|{a}_{q}(\tau)\|_{L^\infty}.
\end{equation}
Since the flow is an homeomorphism then we get again in view of Lemma \ref{mohima1}
\begin{eqnarray}
\label{T52}
\nonumber \|e^{(t-\tau)\Delta}\Delta_{j}\bar{g}_{q}(\tau)\|_{L^\infty}&\lesssim&
  e^{-c(t-\tau)2^{2j}}   \Big(   \Vert [\Delta_{q},v.\nabla]a(\tau)\Vert_{L^{\infty}}\\
  &+&\|(S_{q-1}v-v)\cdot\nabla a_{q}\|_{L^\infty}
\Big).
  \end{eqnarray}
  From Proposition \ref{lemm2} we have
  \begin{eqnarray}\label{TT53}
  \nonumber\Vert [\Delta_{q},v\cdot\nabla]a(t)\Vert_{L^{\infty}}&\lesssim& \|a(t)\|_{L^\infty}\Big(\|\nabla\Delta_{-1}v(t)\|_{L^\infty}+(q+2)\|\omega(t)\|_{L^\infty} \Big) \\
  &\lesssim& \|a^0\|_{L^\infty}\Big(\|\nabla\Delta_{-1}v(t)\|_{L^\infty}+(q+2)\|\omega(t)\|_{L^\infty} \Big).
  \end{eqnarray}
  We have used in the last line the maximum principle: $\|a(t)\|_{L^\infty}\leq\|a^0\|_{L^\infty}$.
  
  On the other hand since $q\in\NN^*,$ we can easily obtain
  \begin{eqnarray}\label{TT54}
  \nonumber\|(S_{q-1}v-v)\cdot\nabla a_{q}\|_{L^\infty}&\lesssim&\|a_{q}\|_{L^\infty} 2^q\sum_{j\geq q-1}2^{-j}\|\Delta_{j}\omega\|_{L^\infty} \\
  &\lesssim&\|a^0\|_{L^\infty}\|\omega\|_{L^\infty}.
  \end{eqnarray}
   Putting together (\ref{ali}), (\ref{Eq:4}), (\ref{T52}), (\ref{TT53}) and (\ref{TT54}) we find
  \begin{eqnarray*}
\|\Delta_{j}\bar{a}_{q}(t)\|_{L^\infty}&\lesssim& e^{-c t 2^{2j}}\|\Delta_{j} a_{q}^0\|_{L^\infty}\\
 &+&
V_{q}(t)e^{CV_{q}(t)} 2^{2q} \int_{0}^te^{-c(t-\tau)2^{2j}}
\|{a}_{q}(\tau)\|_{L^\infty}d\tau\\
&+&(q+2)\|a^0\|_{L^\infty}\int_{0}^te^{-c(t-\tau)2^{2j}}\|\omega(\tau)\|_{L^\infty}d\tau\\
&+&\|a^0\|_{L^\infty}\int_{0}^te^{-c(t-\tau)2^{2j}}\|\nabla\Delta_{-1}v(\tau)\|_{L^\infty}d\tau.
 \end{eqnarray*}
Integrating in time and using  Young inequalities, we obtain for all $j\in\NN$
\begin{eqnarray*}
\|\Delta_{j}\bar{a}_{q}\|_{L^1_{t}L^\infty}&\lesssim &( 2^{2j})^{-1}
\Big(  \|\Delta_{j} a_{q}^0\|_{L^\infty}+(q+2)\|a^0\|_{L^\infty}\|\omega\|_{L^1_{t}L^\infty}+\\
&& \|a^0\|_{L^\infty}\|\nabla\Delta_{-1}v\|_{L^1_{t}L^\infty}\Big)
+V_{q}(t)e^{CV_{q}(t)}2^{2(q-j)}
\|{a}_{q}\|_{L^1_{t}L^\infty}.
\end{eqnarray*}
Let $N$   be a  large integer that will  be chosen later.   Since the flow is an homeomorphism, then we can write
\begin{eqnarray*}
 2^{2q}\|a_{q}\|_{L^1_{t}L^\infty}&=&
 2^{2q}
\|\bar{a}_{q}\|_{L^1_{t}L^\infty}\\
&\leq& 2^{2q}\Big(\sum_{|j-q|<N}
\|\Delta_{j}\bar{a}_{q}\|_{L^1_{t}L^\infty}+
\sum_{|j-q|\geq N}\|\Delta_{j}\bar{a}_{q}\|_{L^1_{t}L^\infty}\Big).
\end{eqnarray*}
Hence, for  all $q> N$, one has
\begin{eqnarray*}
 2^{2q}\|a_{q}\|_{L^1_{t}L^\infty}&\lesssim&
\| a^0\|_{L^\infty}+2^{2N}\|a^0\|_{L^\infty}\Big( (q+2)\|\omega
\|_{L^1_{t}L^\infty}+\|\nabla\Delta_{-1}v\|_{L^1_{t}L^\infty}\Big)\\
&+&V_{q}(t)e^{CV_{q}(t)}2^{2N}
2^{2q}\|a_{q}\|_{L^1_{t}L^\infty}+2^{2q}\sum_{|j-q|\geq N}\|\Delta_{j}\bar{a}_{q}\|_{L^1_{t}L^\infty}.
\end{eqnarray*}
According to Lemma \ref{Vishik00}, we have
$$
\|\Delta_{j}\bar{a}_{q}(t)\|_{L^\infty}\lesssim 2^{-|q-j|}e^{CV_{q}(t)}\|a_{q}(t)\|_{L^\infty}.
$$
Thus, we infer
\begin{eqnarray*}
 2^{2q}\|a_{q}\|_{L^1_{t}L^\infty}&\lesssim&
\| a^0\|_{L^\infty}+2^{2N}\|a^0\|_{L^\infty}\Big( (q+2)\|\omega
\|_{L^1_{t}L^\infty}+\|\nabla\Delta_{-1}v\|_{L^1_{t}L^\infty}\Big)\\
&+&V_{q}(t)e^{CV_{q}(t)}2^{2N}
2^{2q}\|a_{q}\|_{L^1_{t}L^\infty}+2^{-N}e^{CV_{q}(t)}2^{2q}\|{a}_{q}\|_{L^1_{t}L^\infty}.
\end{eqnarray*}

For low frequencies, $q\leq N,$ we write 
$$
 2^{2q}\|a_{q}\|_{L^1_{t}L^\infty}\lesssim
 2^{2N}
\|a\|_{L^1_{t}L^\infty}.
$$
Therefore we get for $q\in\NN\cup\{-1\},$

\begin{eqnarray*}
 2^{2q}\|a_{q}\|_{L^1_{t}L^\infty}&\lesssim&\|a^0\|_{L^\infty}+2^{2N}\|a\|_{L^1_{t}L^\infty}\\
&+&2^{2N}\|a^0\|_{L^\infty}\Big( (q+2)\|\omega
\|_{L^1_{t}L^\infty}+\|\nabla\Delta_{-1}v\|_{L^1_{t}L^\infty}\Big)\\
&+&\Big(V_{q}(t)e^{CV_{q}(t)}2^{2N}
+2^{-N}e^{CV_{q}(t)}\Big)2^{2q}\|{a}_{q}\|_{L^1_{t}L^\infty}.
\end{eqnarray*}
Choosing $N$ and $t$ such that
$$
V_{q}(t)e^{CV_{q}(t)}2^{2N}+e^{CV_{q}(t)} 2^{-N}\lesssim \epsilon,
$$
where $\epsilon<<1.$ This is possible for small time $t$ such that 
$$
V_{q}(t)\leq C_{1},
$$
where $C_{1}$ is a small  absolute constant.

 Under this assumption, one obtains for $q\geq -1$

\begin{eqnarray*}
 2^{2q}\|a_{q}\|_{L^1_{t}L^\infty}&\lesssim&\|a\|_{L^1_{t}L^\infty}+
\| a^0\|_{L^\infty}\Big(1+ (q+2)\|\omega
\|_{L^1_{t}L^\infty}+\|\nabla\Delta_{-1}v\|_{L^1_{t}L^\infty}\Big).
\end{eqnarray*}

Let us now see how to extend this for arbitrarly large time $T.$ We take a partition $(T_{i})_{i=1}^M$ of $[0,T]$ such that
$$
\int_{T_{i}}^{T_{i+1}}\|\nabla S_{q-1}v(t)\|_{L^\infty}dt\simeq C_{1}.
$$
 Reproducing the same arguments as above we find in view of $\|a(T_{i})\|_{L^\infty}\leq \|a^0\|_{L^\infty},$

\begin{eqnarray*}
 2^{2q}\int_{T_{i}}^{T_{i+1}}\|a_{q}(t)\|_{L^\infty}dt&\lesssim&\int_{T_{i}}^{T_{i+1}}\|a(t)\|_{L^\infty}dt+
\| a^0\|_{L^\infty}\\
&+&\|a^0\|_{L^\infty}\Big( (q+2)\int_{T_{i}}^{T_{i+1}}\|\omega(t)
\|_{L^\infty}dt\\
&+&\int_{T_{i}}^{T_{i+1}}\|\nabla\Delta_{-1}v(t)\|_{L^\infty}dt\Big).
\end{eqnarray*}
 Summing these estimates we get
 \begin{eqnarray*}
 2^{2q}\|a_{q}\|_{L^1_{T}L^\infty}&\lesssim&\|a\|_{L^1_{T}L^\infty}+
(M+1)\| a^0\|_{L^\infty}+\\
&+&\|a^0\|_{L^\infty}\Big( (q+2)\|\omega
\|_{L^1_{T}L^\infty}+\|\nabla\Delta_{-1}v\|_{L^1_{T}L^\infty}\Big).
\end{eqnarray*} 
   As $M\approx V_{q}(T),$ then
   \begin{eqnarray*}
 2^{2q}\|a_{q}\|_{L^1_{T}L^\infty}&\lesssim&\|a\|_{L^1_{T}L^\infty}+
(V_{q}(T)+1)\| a^0\|_{L^\infty}+\\
&+&\|a^0\|_{L^\infty}\Big( (q+2)\|\omega
\|_{L^1_{T}L^\infty}+\|\nabla\Delta_{-1}v\|_{L^1_{T}L^\infty}\Big).
\end{eqnarray*} 
Since 
$$
\|\nabla S_{q-1} v\|_{L^\infty}\leq \|\nabla\Delta_{-1}v\|_{L^\infty}+(q+2)\|\omega\|_{L^\infty},
$$
then inserting this estimate into the previous one
\begin{eqnarray*}
 2^{2q}\|a_{q}\|_{L^1_{T}L^\infty}&\lesssim&\|a^0\|_{L^\infty}\Big((1+T)+ (q+2)\|\omega
\|_{L^1_{T}L^\infty}+\|\nabla\Delta_{-1}v\|_{L^1_{T}L^\infty}\Big).
\end{eqnarray*} 
This is the desired result.
 \end{proof}
 \section{Proof of Theorem \ref{thmm0}}
 We restrict ourselves to the {\it a priori} estimates. The existence and uniqueness parts are easily obtained  with small modifications of the proof of Theorem \ref{thm1}.
 \begin{prop}\label{prop992}
 For $v^0\in B_{p,1}^{1+\frac 2p}$ and $\theta^0\in L^r,$ with $2<r\leq p,$ we have for $t\in \RR_{+}$
 \begin{itemize}
\item[{\rm 1)}]
  $$
 \|\theta(t)\|_{L^r}\leq\|\theta^0\|_{L^r}.
$$ 
 \item[{\rm 2)}]
 $$
\|\theta\|_{L^1_{t}B_{r,1}^{1+\frac2r}}+\|\omega(¤t)\|_{L^\infty}+\|\omega\|_{\widetilde L^\infty_{t}B_{\infty,1}^0}\le C_{0}e^{e^{C_{0}t}}.
 $$
\item[{\rm 3)}]
  $$\|\theta\|_{\widetilde L^1_{t}B_{r,\infty}^2}+\|v\|_{\widetilde L^\infty_{t }B_{p,1}^{1+\frac 2p}}\le C_{0}e^{e^{e^{C_{0}t}}},
$$
\end{itemize}
where the constant  $C_{0}$ depends on the quantity $\|\theta^0\|_{L^r}$ and $\|v^0\|_{B_{p,1}^{1+\frac2p}}.$
 \end{prop}
 \begin{proof}
The first estimate  can be easily obtained from $L^r$ energy estimate. for the second one, we recall that the  vorticity $\omega=\partial_{1}v^2-\partial_{2}v^1$ satisifies the equation
 \begin{equation}\label{TY}
 \partial_{t}\omega+v\cdot\nabla\omega=\partial_{1}\theta.
 \end{equation}
 Taking the $L^\infty$ norm we get
 \begin{equation}\label{y34}
 \|\omega(t)\|_{L^\infty}\le \|\omega^0\|_{L^\infty}+\|\nabla\theta\|_{L^1_{t}L^\infty}.
 \end{equation}
Using the embedding $B_{r,1}^{1+\frac2r}\hookrightarrow \textnormal{Lip}(\RR^2)$ we obtain \begin{eqnarray*}
 \|\omega(t)\|_{L^\infty}&\lesssim &\|\omega^0\|_{L^\infty}+\|\theta\|_{L^1_{t}B_{r,1}^{1+\frac 2r}}.
\end{eqnarray*}
From Theorem \ref{effect} applied to the temperature equation and by Bernstein inequalities we deduce for $\epsilon>0,$ \begin{eqnarray*}
\|\theta\|_{L^1_{t}B_{r,\infty}^{2-\epsilon}}&\lesssim &\|\theta^0\|_{L^r}\Big(1+t+\|\omega\|_{L^1_{t}L^\infty}+\|\Delta_{-1}\nabla v\|_{L^1_{t}L^\infty}        \Big)\\
&\lesssim&\|\theta^0\|_{L^r}\Big(1+t+\|\omega\|_{L^1_{t}L^\infty}+\|\nabla v\|_{L^1_{t}L^p}        \Big).
\end{eqnarray*}
This leads for $r>2$ to the inequality
$$\|\theta\|_{L^1_{t}B_{r,\infty}^{1+\frac2r}}\lesssim\|\theta^0\|_{L^r}\Big(1+t+\|\omega\|_{L^1_{t}L^\infty}+\|\nabla v\|_{L^1_{t}L^p}        \Big).
$$
On the other hand   we have the classical result $\|\nabla v\|_{L^p}\approx\|\omega\|_{L^p},$ for $p\in]1,\infty[.$ Thus we get
$$
\|\theta\|_{L^1_{t}B_{r,1}^{1+\frac 2r}}\lesssim \|\theta^0\|_{L^r}\Big(1+t+\|\omega\|_{L^1_{t}L^\infty}+\|\omega\|_{L^1_{t}L^p}\Big).
$$
The estimate of the $L^p$ norm of the vorticity can be done as its $L^\infty$ norm ($r\leq p$)
$$
\|\omega(t)\|_{L^p}\lesssim\|\omega^0\|_{L^p}+\|\theta\|_{L^1_{t}B_{r,1}^{1+\frac 2r}}.
$$
Set $f(t):=\|\omega(t)\|_{L^\infty\cap L^p}+\|\theta\|_{L^1_{t}B_{r,1}^{1+\frac 2r}}.$ Then combining the above estimates yields
$$
f(t)\lesssim\|\omega^0\|_{L^\infty\cap L^p}+\|\theta^0\|_{L^r}(1+t)+\|\theta^0\|_{L^r}\int_{0}^tf(\tau)d\tau.
$$
According to Gronwall's inequality, one has
\begin{equation}\label{ZY3}
\|\omega(t)\|_{L^\infty\cap L^p}+\|\theta\|_{L^1_{t}B_{r,1}^{1+\frac 2r}}\lesssim
\Big(\|\omega^0\|_{L^\infty\cap L^p}+\|\theta^0\|_{L^r}(1+t)\Big) e^{Ct\|\theta^0\|_{L^r}}\leq C_{0}e^{C_{0}t},
\end{equation}
where $C_{0}$ is a constant depending on the initial data.

Let us now turn to the estimate of $\|\omega(t)\|_{B_{\infty,1}^{0}}.$ From Proposition \ref{prop11} and Besov embeddings,
\begin{eqnarray}\label{ZY5}
\nonumber\|\omega\|_{\widetilde L^\infty_{t}B_{\infty,1}^0}&\lesssim&\big(\|\omega^0\|_{L^\infty}+\|\theta\|_{L^1_{t}B_{\infty,1}^1}\big)\big(1+\|\nabla v\|_{L^1_{t}L^\infty} \big)\\
&\lesssim&\big(\|\omega^0\|_{L^\infty}+\|\theta\|_{L^1_{t}B_{r,1}^{1+\frac2r}}\big)\big(1+\|\nabla v\|_{L^1_{t}L^\infty} \big).
\end{eqnarray}
On the other hand we have
\begin{eqnarray}\label{ZY4}
\nonumber\|\nabla v(t)\|_{L^\infty}&\leq&\|\nabla\Delta_{-1}v(t)\|_{L^\infty}+\sum_{q\in\NN}\|\Delta_{q}\nabla v(t)\|_{L^\infty}\\
\nonumber&\lesssim&\|\nabla\Delta_{-1}v(t)\|_{L^p}+\|\omega(t)\|_{B_{\infty,1}^0}\\
&\lesssim&\|\omega(t)\|_{L^p}+\|\omega\|_{\widetilde L^\infty_{t}B_{\infty,1}^0}.
\end{eqnarray}
Putting together (\ref{ZY3}), (\ref{ZY5}) and (\ref{ZY4}) and using Gronwall's inequality gives
\begin{equation}\label{ZY7}
\|\nabla v\|_{L^\infty}+\|\omega\|_{\widetilde L^\infty_{t}B_{\infty,1}^0}\le C_{0}e^{e^{C_{0}t}}.
\end{equation}
It remains to prove the third point of the proposition. The smoothing effect on  $\theta$ is a direct consequence of  (\ref{visc0}) and the above inequality,
$$
\|\theta\|_{\widetilde L^1_{t}B_{r,\infty}^2}\le C_{0}e^{e^{C_{0}t}}.
$$
Concerning the velocity estimate we write
$$
\|v\|_{\widetilde L^\infty_{t}B_{p,1}^{1+\frac 2p}}\lesssim \|v\|_{L^\infty_{t}L^p}+\|\omega\|_{\widetilde L^\infty_{t}B_{p,1}^{\frac 2p}}.
$$
Using the velocity equation, we obatin
$$
\|v(t)\|_{L^p}\le \|v^0\|_{L^p}+t\|\theta^0\|_{L^p}+\int_{0}^t\|\mathcal{P}(v\cdot\nabla v)(\tau)\|_{L^p}d\tau.
$$
where $\mathcal{P}$ denotes Leray projector. It follows from classical estimate that
$$
\|\mathcal{P}(v\cdot\nabla v)\|_{L^p}\lesssim\|v\cdot\nabla v\|_{L^p}\lesssim\|v\|_{L^p}\|\nabla v\|_{L^\infty}.
$$
Thus we get in view of Gronwall's inequality and (\ref{ZY7})
\begin{equation}\label{ZY8}
\|v\|_{L^\infty_{t}L^p}\leq C_{0}e^{e^{e^{C_{0}t}}}.
\end{equation}
It remains to estimate $\|\omega(t)\|_{B_{p,1}^{\frac 2p}}.$ We apply Proposition \ref{propagation} to the vorticity equation and we use Besov embeddings,
\begin{eqnarray*}
\|\omega\|_{\widetilde L^\infty_{t}B_{p,1}^{\frac2p}}&\lesssim& e^{CV(t)}(\|\omega^0\|_{B_{p,1}^{\frac2p}}+\|\theta\|_{L^1_{t}B_{p,1}^{1+\frac2p}})
\\
&\lesssim&e^{CV(t)}(\|\omega^0\|_{B_{p,1}^{\frac2p}}+\|\theta\|_{L^1_{t}B_{r,1}^{1+\frac2r}}).
\end{eqnarray*}
It suffices now to use (\ref{ZY3}) and (\ref{ZY7}).
\end{proof}
 \section{Proof of Theorem \ref{thm1}}
 The case $p=+\infty$ is more subtle and the difficulty comes from the term $\|\nabla\Delta_{-1}v\|_{L^\infty},$ since Riesz transforms do not map $L^\infty$ to itself. To avoid this problem we use a frequency interpolation method. The proof will be done in several steps. The first one deals with  some {a priori} estimates. We give in the second the uniqueness result and the last is reserved to the proof of the existence part. 
 \subsection{A priori estimates}
 The main result of this section is the following:
 \begin{prop}\label{prop1000}
 There exists a constant $C_{0}$ depending on $\|v^0\|_{B_{\infty,1}^0}$ and $\|\theta^0\|_{L^\infty}$ such that for $t\in[0,\infty[$
 \begin{eqnarray*}
\|\theta(t)\|_{L^\infty}\leq\|\theta^0\|_{L^\infty};\quad \|\theta\|_{L^1_{t}B_{\infty,1}^1}&\leq& C_{0}e^{C_{0}t^3}\quad\hbox{and}\\
\|v\|_{\widetilde L^\infty_{t}B_{\infty,1}^1}+\|\theta\|_{\widetilde L^1_{t}B_{\infty,\infty}^2}+\|\theta(t)\|_{\mathcal{B}^\infty}&\leq &C_{0}e^{e^{C_{0}t^3}}.
\end{eqnarray*}
 \end{prop}
 \begin{proof}
 The $L^\infty$-bound of the temperature can be  easily obtained from the maximum principle. To give the other bounds we start with the following estimate for the vorticity, which is again a direct consequence of the maximum principle,
\begin{equation}\label{Tik}
\|\omega(t)\|_{L^\infty}\leq\|\omega^0\|_{L^\infty}+\|\nabla\theta\|_{L^1_{t}L^\infty}\leq\|\omega^0\|_{L^\infty}+\|\theta\|_{L^1_{t}B_{\infty,1}^1}.
\end{equation}
Let $N\in\NN^*,$ then we get by definition of Besov spaces and the maximum principle
\begin{eqnarray*}
\|\theta\|_{L^1_{t}B_{\infty,1}^1}&=& \sum_{q\leq N-1}2^q\|\Delta_{q}\theta\|_{L^1_{t}L^\infty}+\sum_{q\geq N}2^q\|\Delta_{q}\theta\|_{L^1_{t}L^\infty}\\
&\lesssim& 2^Nt\|\theta^0\|_{L^\infty}+\sum_{q\geq N}2^q\|\Delta_{q}\theta\|_{L^1_{t}L^\infty}.
\end{eqnarray*}
By virtue of Theorem {\ref{effect}} one has
\begin{eqnarray*}
\nonumber\|\theta\|_{L^1_{t}B_{\infty,1}^1}&\lesssim& 2^{N}t\|\theta^0\|_{L^\infty}+2^{-N}\|\theta^0\|_{L^\infty}\Big( 1+t+\|\nabla\Delta_{-1}v\|_{L^1_{t}L^\infty}
+N\|\omega\|_{L^1_{t}L^\infty}\Big)
\\
&\lesssim&2^{N}t\|\theta^0\|_{L^\infty}+2^{-N}\|\theta^0\|_{L^\infty}\Big( 1+t+\|\nabla\Delta_{-1}v\|_{L^1_{t}L^\infty}\Big)+\|\omega\|_{L^1_{t}L^\infty}.
\end{eqnarray*}
Choosing judiciously $N$ we get
\begin{equation}
\label{tout01}
\|\theta\|_{L^1_{t}B_{\infty,1}^1}\lesssim \|\omega\|_{L^1_{t}L^\infty}+t^{\frac{1}{2}}\|\theta^0\|_{L^\infty}\Big(1+t+\|\nabla\Delta_{-1}v\|_{L^1_{t}L^\infty}\Big)^{\frac{1}{2}}
\end{equation}
The following lemma gives an estimate of  the low frequency of the velocity.
\begin{lem}\label{LM1}
For all $t\geq 0,$ wed have
$$
\|\nabla\Delta_{-1}v(t)\|_{L^\infty}\lesssim1+\log\big(e+\|v^0\|_{L^\infty}+t\|\theta^0\|_{L^\infty}  \big)\|\omega\|_{L^\infty_{t}L^\infty}+t\|\omega\|_{L^\infty_{t}L^\infty}^2.
$$
\end{lem}
\begin{proof}
 Fix $N\in\NN^*.$ Since $\Delta_{-1}=\Delta_{-1}(\dot{S}_{-N}+\sum_{q=-N}^0\dot\Delta_{q})$ then we have
\begin{eqnarray*}
\|\nabla\Delta_{-1}v\|_{L^\infty}&\lesssim&\|\nabla\dot S_{-N}v\|_{L^\infty}+\sum_{q=-N}^0\|\nabla\dot\Delta_{q}v\|_{L^\infty}\\
&\lesssim&2^{-N}\|v\|_{L^\infty}+\sum_{-N}^0\|\dot\Delta_{q}\omega\|_{L^\infty}\\
&\lesssim& 2^{-N}\|v\|_{L^\infty}+N\|\omega\|_{L^\infty}.
\end{eqnarray*}
Taking $N\approx \log(e+\|v\|_{L^\infty})$ we get
\begin{equation}\label{Tout11}
\|\nabla\Delta_{-1}v\|_{L^\infty}\lesssim 1+\|\omega\|_{L^\infty}\log(e+\|v\|_{L^\infty}).
\end{equation}
It remains to estimate $\|v\|_{L^\infty}.$ Let $M\in\NN$ then we have
$$
\|v\|_{L^\infty}\lesssim \|\dot S_{-M}v\|_{L^\infty}+2^M\|\omega\|_{L^\infty}.
$$
Now using the equation of the velocity we get

\begin{eqnarray*}
\|\dot S_{-M}v(t)\|_{L^\infty}&\leq&\|\dot S_{-M}v^0\|_{L^\infty}+\|\dot S_{-M}\theta\|_{L^1_{t}L^\infty}\\
&+&\int_{0}^t\|\dot S_{-M}\diver\mathcal{P}(v\otimes v)(\tau)\|_{L^\infty}d\tau\\
&\lesssim&\|v^0\|_{L^\infty}+t\|\theta^0\|_{L^\infty}+2^{-M}\int_{0}^t\|v(\tau)\|_{L^\infty}^2d\tau.
\end{eqnarray*}
We have used the following inequality
$$
\|\dot S_{-M}\diver\mathcal{P}(v\otimes v)\|_{L^\infty}\leq\sum_{q\leq -M-1}\|\dot \Delta_{q}\diver\mathcal{P}(v\otimes v)\|_{L^\infty}\lesssim \sum_{q\leq -M-1}2^q\|v\otimes v\|_{L^\infty}.
$$
Thus we obtain
$$
\|v\|_{L^\infty}\lesssim \|v^0\|_{L^\infty}+t\|\theta^0\|_{L^\infty}+2^{-M}\int_{0}^t\|v(\tau)\|_{L^\infty}^2d\tau+2^M\|\omega(t)\|_{L^\infty}.
$$
Taking $M$ such that 
$$
2^{2M}\approx\frac{\int_{0}^t\|v\|_{L^\infty}^2d\tau}{\|\omega\|_{L^\infty}},
$$
we find
$$\|v\|_{L^\infty}\lesssim \|v^0\|_{L^\infty}+t\|\theta^0\|_{L^\infty}+\|\omega(t)\|_{L^\infty}^{\frac{1}{2}}\Big(\int_{0}^t\|v(\tau)\|_{L^\infty}^2d\tau\Big)^{\frac{1}{2}}.
$$
According to Gronwall's inequality we get
\begin{equation}\label{tr67}\|v\|_{L^\infty}\lesssim \big(\|v^0\|_{L^\infty}+t\|\theta^0\|_{L^\infty}\big)e^{Ct\|\omega\|_{L^\infty_{t}L^\infty}}.
\end{equation}
Inserting this estimate into (\ref{Tout11}) we find the desired inequality.
\end{proof}
 Lemma \ref{LM1} and (\ref{tout01}) yield
\begin{eqnarray*}
\nonumber\|\theta\|_{L^1_{t}B_{\infty,1}^1}^2&\leq&C_{0}(1+t^2)+\|\omega\|_{L^1_{t}L^\infty}^2+C_{0}(1+t^2)\int_{0}^t\|\omega\|_{L^\infty_{\tau}L^\infty}^2d\tau\\
&\leq&C_{0}(1+t^2)\Big( 1+\int_{0}^t\|\omega\|_{L^\infty_{\tau}L^\infty}^2d\tau\Big).
\end{eqnarray*}
Combining this estimate with (\ref{Tik}) yields
$$
\|\omega\|_{L^\infty_{t}L^\infty}^2\leq C_{0}(1+t^2)\Big( 1+\int_{0}^t\|\omega\|_{L^\infty_{\tau}L^\infty}^2d\tau\Big).
$$
Applying Gronwall's inequality we get
\begin{equation}\label{Tik1}
\|\omega(t)\|_{L^\infty}\leq C_{0}e^{C_{0}t^3}.
\end{equation}
This gives
\begin{equation}\label{Tik2}
\|\theta\|_{L^1_{t}B_{\infty,1}^1}\leq C_{0}e^{C_{0}t^3}.
\end{equation}
From Lemma \ref{LM1} we have
\begin{equation}\label{Tik55}
\|\nabla\Delta_{-1}v(t)\|_{L^\infty}\leq C_{0}e^{C_{0}t^3}.
\end{equation}

\

Let us now turn to the estimate of the vorticity in $B_{\infty,1}^0$ space. For this purpose we apply Proposition \ref{prop11} to the vorticity equation, with $p=+\infty$ and $r=1$
\begin{equation}\label{Tik56}
\|\omega\|_{\widetilde L^\infty_{t}B_{\infty,1}^0}\lesssim\big(\|\omega^0\|_{L^\infty}+\|\nabla\theta\|_{L^1_{t}B_{\infty,1}^0}\big)\Big(1+\int_{0}^t\|\nabla v(\tau)\|_{L^\infty}d\tau\Big).
\end{equation}
On the other hand we have by definition and from (\ref{Tik55}) and (\ref{Tik56})
\begin{eqnarray*}
\| \nabla v(t)\|_{L^\infty}\lesssim\|v\|_{\widetilde L^\infty_{t}B_{\infty,1}^1}&\lesssim& \|\nabla\Delta_{-1}v\|_{L^\infty_{t}L^\infty}+\sum_{q\in\NN}\|\Delta_{q}\omega\|_{L^\infty_{t}L^\infty}\\
&\lesssim&C_{0}e^{C_{0}t^3}+\|\omega\|_{\widetilde L^\infty_{t}B_{\infty,1}^0}\\
&\lesssim& C_{0}e^{C_{0}t^3}\Big(1+\int_{0}^t\| v(\tau)\|_{B_{\infty,1}^1}d\tau   \Big).
\end{eqnarray*}
It suffices now  to use Gronwall's inequality.

\

To estimate $\|\theta\|_{L^1_{t}B_{\infty,\infty}^2}$ it suffices to combine (\ref{visc0}) with the Lipschitz estimate of the velocity. The last estimate $\|\theta(t)\|_{\mathcal{B}^\infty}$ will be done as follows:
$$
\|\theta(t)\|_{\mathcal{B}^\infty}\leq\|\theta^0\|_{L^\infty}+\sum_{q\leq 0}\|\dot{\Delta}_{q}\theta(t)\|_{L\infty}.
$$
Using the temperature equation we find 
\begin{eqnarray*}
\|\dot\Delta_{q}\theta(t)\|_{L^\infty}&\leq&\|\dot\Delta_{q}\theta^0\|_{L^\infty}+\|\dot{\Delta}_{q}(v\cdot\nabla\theta)\|_{L^1_{t}L^\infty}+\|\dot\Delta_{q}\Delta\theta\|_{L^1_{t}L^\infty}\\
&\lesssim&\|\dot\Delta_{q}\theta^0\|_{L^\infty}+2^q\|v\,\theta\|_{L^1_{t}L^\infty}+2^{2q}\|\theta\|_{L^1_{t}L^\infty}\\
&\lesssim&\|\dot\Delta_{q}\theta^0\|_{L^\infty}+2^q\|\theta^0\|_{L^\infty}\|v\|_{L^1_{t}L^\infty}+2^{2q}t\|\theta^0\|_{L^\infty}.
\end{eqnarray*}
Therefore we get 
$$
\sum_{q\leq 0}\|\dot{\Delta}_{q}\theta(t)\|_{L\infty}\lesssim\sum_{q\leq 0}\|\dot{\Delta}_{q}\theta^0\|_{L\infty}+C_{0}e^{e^{C_{0}t^3}}.
$$
This concludes the proof of the proposition.
\end{proof}
\subsection{Uniqueness part}
As it was shown in the previous paragraph we can give  an {\it a priori} estimates for both Lipschitz norms of the velocity and the temperature only under the assumption $v^0\in B_{\infty,1}^0$ and $\theta^0\in L^\infty.$ However it seems that  the uniqueness part (even the existence) needs an addition condition of the initial data $\theta^0,$ \mbox{namely $\Delta_{-1}\in \dot B_{\infty,1}^0.$ }

                 Let us consider two solutions $\{(v^{j},\theta^{j})\}_{j=1}^2$
         for the system $(B_{0,\kappa}),$  with initial data $(v^{j,0},\theta^{j,0}),\, j=1,2$ and     satisfying for a fixed time $T>0$
         $$v^{j}\in L^\infty_{T}\,B_{\infty,1}^1\quad\quad\hbox{and}\quad\quad\theta^{j}\in L^\infty_{T}\mathcal{B}^{\infty}\cap L^1_{T}\,\textnormal{Lip}(\RR^2).
         $$
          We set
         $$
         v=v^{1}-v^{2},\, \theta=\theta^{1}-\theta^{2},\,\pi=\pi^1-\pi^2,\, v^0=v^{1,0}-v^{2,0},
         \quad\hbox{and}\quad\theta^0=\theta^{1,0}-\theta^{2,0}.$$
        Thus we have the equations
        \begin{equation}\label{eqs23}
        \partial_{t} v+v^{1}\cdot\nabla v=-\nabla\pi- v\cdot\nabla v^{2}+\theta e_{2},
        \end{equation}   
        \begin{equation}\label{eqs24}     \partial_{t}\theta+v^{1}\cdot\nabla\theta-\Delta\theta=-v\cdot\nabla \theta^{2}.
        \end{equation}
    Hereafter we denote $V_{j}(t):=C\|\nabla v^j\|_{L^1_{t}L^\infty}, j=1,2.$ Now applying 
 Lemma \ref{propagation} yields     $$
      \|v(t)\|_{B_{\infty,1}^0}\lesssim e^{V_{1}(t)}\Big(\|v^0\|_{B_{\infty,1}^0}+\int_{0}^t\big(\|\nabla\pi(\tau)\|_{B_{\infty,1}^0}+\|v\cdot\nabla v^2(\tau)\|_{B_{\infty,1}^0}+      \|\theta(\tau)\|_{B_{\infty,1}^0}\big)d\tau
      \Big).    $$   
      To estimate the pressure we write the following identity
      $$
      \Delta\pi=-\textnormal{div }\big( v^1\cdot\nabla v+v\cdot\nabla v^2 \big)+\partial_{2}\theta.
      $$
      Since $\textnormal{div }\big( v^1\cdot\nabla v \big)=\textnormal{div }\big( v\cdot\nabla v^1 \big)$ then 
     $$ \nabla\pi=-\nabla\Delta^{-1}\textnormal{div }\big( v\cdot\nabla(v^1+v^2)  \big)+\nabla\Delta^{-1}\partial_{2}\theta.
     $$     
      From the embedding $\dot B_{\infty,1}^0\hookrightarrow B_{\infty,1}^0$ and the fact that Riesz transforms act continuously on homogeneous Besov spaces, one obtains
      $$
      \|\nabla\pi\|_{B_{\infty,1}^0}\lesssim \|v\cdot\nabla(v^1+v^2)\|_{\dot B_{\infty,1}^0}+\|\theta\|_{\dot B_{\infty,1}^0}.
      $$
      It is easy to see that
      \begin{eqnarray*}
      \|v\cdot\nabla(v^1+v^2)\|_{\dot B_{\infty,1}^0}&\lesssim &\sum_{q\leq0}2^q\|\dot\Delta_{q}(v\otimes(v^1+v^2))\|_{L^\infty}+\sum_{q>0}\|\Delta_{q}(v\cdot\nabla(v^1+v^2))\|_{L^\infty}\\
      &\lesssim&
      \|v\|_{L^\infty}\|v^1+v^2\|_{L^\infty}+\|v\cdot\nabla(v^1+v^2)\|_{ B_{\infty,1}^0}.
      \end{eqnarray*}
      Using Bony's decomposition and the incompressibilty of the velocity $v$ one obtains the general estimate
      $$
      \|v\cdot\nabla w\|_{B_{\infty,1}^0}\lesssim \|v\|_{B_{\infty,1}^0}\|w\|_{B_{\infty,1}^1}.
      $$
      It follows that
      \begin{eqnarray*}
      \|v\cdot\nabla(v^1+v^2)\|_{ B_{\infty,1}^0}&\lesssim&\|v\|_{B_{\infty,1}^0}\|v^1+v^2\|_{B_{\infty,1}^1}\\
    \|v\cdot\nabla v^2\|_{ B_{\infty,1}^0}&\lesssim& \|v\|_{B_{\infty,1}^0}\|v^2\|_{B_{\infty,1}^1}.
      \end{eqnarray*}
     Putting together these estimates gives
         \begin{eqnarray}\label{ham77}
      \|v(t)\|_{B_{\infty,1}^0}&\lesssim& e^{V_{1}(t)}\Big(\|v^0\|_{B_{\infty,1}^0}+ \|\theta\|_{L^1_{t}\dot B_{\infty,1}^0} +   \int_{0}^t\|v(\tau)\|_{B_{\infty,1}^0}w_{1,2}(\tau)d\tau  \Big),
        \end{eqnarray}
     with   $$
        w_{1,2}(t):=\|v^1(t)\|_{B_{\infty,1}^1}  +\|v^2(t)\|_{B_{\infty,1}^1}.
               $$
       It remains to estimate the quantity $\|\theta\|_{L^1_{t}\dot B_{\infty,1}^0}.$ For this aim we split $\theta$  into low and high frequencies 
       $$
       \|\theta\|_{\dot B_{\infty,1}^0}\leq\sum_{q\leq 0}\|\dot\Delta_{q}\theta\|_{L^\infty}+\|\theta-\Delta_{-1}\theta\|_{B_{\infty,1}^0}.
       $$
       Using the equation of $\theta$ we get easily
       \begin{eqnarray*} \|\dot\Delta_{q}\theta(t)\|_{L^\infty}&\leq&\|\dot\Delta_{q}\theta^{0}\|_{L^\infty}+\|\dot\Delta_{q}(\textnormal{div\,}(v^1\;\theta+v\;\theta^2)\|_{L^1_{t}L^\infty}+\|\dot\Delta_{q}\Delta\theta\|_{L^1_{t}L^\infty}\\
      & \lesssim&\|\dot\Delta_{q}\theta^{0}\|_{L^\infty}+2^q\int_{0}^t\big( \|v^1(\tau)\|_{L^\infty}\|\theta(\tau)\|_{L^\infty}+\|\theta^2(\tau)\|_{L^\infty}\|v(\tau)\|_{L^\infty}\big)d\tau\\
       &+&2^{2q}\|\theta\|_{L^1_{t}L^\infty}     \end{eqnarray*}
                Combining both last estimates with Besov embeddings (essentially $\dot B_{\infty,1}^0\hookrightarrow B_{\infty,1}^0\hookrightarrow L^\infty$) gives
         \begin{eqnarray*}
         \|\theta(t)\|_{\dot B_{\infty,1}^0}&\lesssim &\sum_{q\leq0}\|\dot\Delta_{q}\theta^0\|_{L^\infty}+ \int_{0}^t(1+w_{1,2}(\tau))\|\theta(\tau)\|_{\dot B_{\infty,1}^0}d\tau\\
         &+&\|\theta^{2,0}\|_{L^\infty}         \int_{0}^t\|v(\tau)\|_{B_{\infty,1}^0}d\tau 
         +\|\theta(t)-\Delta_{-1}\theta(t)\|_{B_{\infty,1}^0}.        
         \end{eqnarray*}
     Integrating over the time and using Gronwall's inequality
     \begin{eqnarray}\label{rt23}
       \nonumber  \|\theta\|_{L^1_{t}\dot B_{\infty,1}^0}\lesssim e^{Ct+Ct\|w_{1,2}\|_{L^\infty_{t}}}\Big(\sum_{q\leq0}\|\dot\Delta_{q}\theta^0\|_{L^\infty}
                  &+&\|\theta^{2,0}\|_{L^\infty}         \int_{0}^t\|v(\tau)\|_{B_{\infty,1}^0}d\tau \\
                   &
         +&\|\theta-\Delta_{-1}\theta\|_{L^1_{t}B_{\infty,1}^0} \Big).       
         \end{eqnarray}        
          It remains to estimate $\|\theta-\Delta_{-1}\theta\|_{L^1_{t}B_{\infty,1}^0}.$ For this purpose we apply Proposition \ref{propagation} to the equation (\ref{eqs24}) with $s=-\frac{1}{2}$ and $\,p=r=\infty$
         \begin{eqnarray*}
         \|\theta-\Delta_{-1}\theta\|_{L^1_{t}{B_{\infty,\infty}^{\frac32}}}&\lesssim& e^{V_{1}(t)}\Big(\|\theta^0\|_{B_{\infty,\infty}^{\frac{-1}{2}}}+\int_{0}^t\|v\cdot\nabla \theta^2(\tau)\|_{L^\infty}d\tau\Big).
               \end{eqnarray*}
    We have used in the above inequality the embedding $L^\infty\hookrightarrow B_{\infty,\infty}^{-\frac12}.$        
    
     Since $B_{\infty,\infty}^{\frac32}\hookrightarrow B_{\infty,1}^0,$ we find
                 \begin{eqnarray*}
         \|\theta-\Delta_{-1}\theta\|_{L^1_{t}{B_{\infty,\infty}^{0}}}&\lesssim& e^{Ct\|w_{1,2}\|_{L^\infty_{t}}}\Big(\|\theta^0\|_{B_{\infty,\infty}^{\frac{-1}{2}}}+\int_{0}^t\|v(\tau)\|_{B_{\infty,1}^0}\|\nabla \theta^2(\tau)\|_{L^\infty}d\tau\Big).
               \end{eqnarray*}         Inserting this estimate into (\ref{rt23}) we get
      
         \begin{eqnarray*}
         \|\theta\|_{L^1_{t}\dot B_{\infty,1}^0}&\lesssim &
e^{Ct+C\|w_{1,2}\|_{L^\infty_{t}}}\Big(\sum_{q\leq0}\|\dot\Delta_{q}\theta^0\|_{L^\infty}+
\|\theta^0\|_{{B}_{\infty,\infty}^{\frac{-1}{2}}}\\
&+&\int_{0}^t\|v(\tau)\|_{B_{\infty,1}^0}\big(\|\nabla\theta^2(\tau)\|_{L^\infty}+\|\theta^{2,0}\|_{L^\infty}\big)d\tau\Big).    
             \end{eqnarray*}  
                       
    Putting together this estimate with (\ref{ham77}) we obtain
    \begin{eqnarray*}
      \|v(t)\|_{B_{\infty,1}^0}&\lesssim&e^{Ct+Ct\|w_{1,2}\|_{L^\infty_{t}}}\Big(\sum_{q\leq 0}\|\dot\Delta_{q}\theta^0\|_{L^\infty}+\|\theta^0\|_{B_{\infty,\infty}^{\frac{-1}{2}}}+\|v^0\|_{B_{\infty,1}^0}\\
      &+&\int_{0}^t\|v(\tau)\|_{B_{\infty,1}^0}\big( w_{1,2}(\tau)+\|\nabla\theta^2(\tau)\|_{L^\infty}+\|\theta^{2,0}\|_{L^\infty}\big)d\tau  \Big).
        \end{eqnarray*}
       It follows from both last estimates and Gronwall's inequality
       \begin{equation}\label{unicite1}
       \|v\|_{L^\infty_{t}B_{\infty,1}^0}+\|\theta\|_{L^1_{t}\dot B_{\infty,1}^0}\leq \eta(t)\Big(\sum_{q\leq 0}\|\dot\Delta_{q}\theta^0\|_{L^\infty}+\|\theta^0\|_{B_{\infty,\infty}^{\frac{-1}{2}}}+\|v^0\|_{B_{\infty,1}^0}\Big),             \end{equation}
       where $\eta=\RR_{+}\to\RR_{+}$  is a function depending on the quantities $\|v^j\|_{L^\infty_{t}B_{\infty,1}^1},$ $\|\nabla\theta^j\|_{L^1_{t}L^\infty}$ and $\|\theta^{j,0}\|_{L^\infty}$.           This concludes the proof of the uniqueness part. 
\subsection{Existence part}
  We will briefly outline the proof of the existence part which is classical. 
    We smooth out  the initial data $(v_{n}^0,\theta_{n}^0):=(S_{n}v^0,S_{n}\theta^0),$ which is uniformly bounded in $B_{\infty,1}^1\times \mathcal{B}^\infty$ and it is easy to check the following convergence result
    $$
   \lim_{n\to\infty} \Big(\|v_{n}^0-v^0\|_{B_{\infty,1}^0}+\sum_{q\leq 0}\|\dot\Delta_{q}(\theta_{n}^0-\theta^0)\|_{L^\infty}+\|\theta^0_{n}-\theta^0\|_{B_{\infty,\infty}^{\frac{-1}{2}}}\Big)=0.    $$
   Since the initial data $(v_{n}^0,\theta^0_{n})$ are smooth then the corresponding Boussinesq system has  global unique smooth solutions  $(v_{n},\theta_{n}).$   In view of  Proposition \ref{prop1000} one has  for every $t\in\RR_{+},$   the uniform estimates:
     $$
\|v_{n}\|_{L^\infty_{t}B_{\infty,1}^1}+\|\theta_{n}(t)\|_{L^\infty}+\|\nabla \theta_{n}\|_{L^1_{t}L^\infty}\leq C_{0} e^{e^{C_{0}t^3}}.    $$
   Now according to (\ref{unicite1}) the sequence $(v_{n},\theta_{n})_{n}$ 
    converges strongly
    in $L^\infty_{\textnormal{loc}}(\RR_{+};B_{\infty,1}^0)\times L^1_{\textnormal loc}(\RR_{+}; \dot B_{\infty,1}^0)$ to $(v,\theta).$ This is sufficiently to pass 
    to the limit in the equations and deduce that $(v,\theta)$ satisfies the system $(B_{0,\kappa}).$\\
    It remains to show  the continuity in time of the velocity. This comes from the estimate $\|v\|_{\widetilde L^\infty_{t}B_{\infty,1}^1},$ (see \cite{T} for more further details).
     
     The proof of Theorem \ref{thm1} is now complete.   
%\appendix
\section*{Appendix A.  Logarithmic estimate}
We shall now give a  logarithmic estimate which is an extension of  Vishik's one \cite{v1}. Our result was firstly proved in \cite{ST} and for the convenience of the reader we will give here the proof.     \begin{prop}\label{prop11}
Let $p,\,r\in[1,+\infty],$ $v$ be a divergence-free vector field belonging to the space $L^1_{loc}(\RR_{+};\Lip(\RR^d))$ and let $a$ be a smooth solution of 
the following equation $($with $\nu\geq 0$$)$,
$$\left\lbrace
\begin{array}{l}
\partial_t a+v\cdot\nabla a- \nu\Delta a=f
\\
{a}_{| t=0}=a^{0}.\\
\end{array}
\right.
$$
If the initial data  $a^{0}\in B_{p,r}^0,$ then we have for all $t\in\RR_{+}$

$$
\|a\|_{\widetilde L^\infty_{t}B_{p,r}^0}\leq C\big(\|a^{0}\|_{B_{p,r}^0}+\|f\|_{\widetilde L^1_{t}B_{p,r}^0}\big)
\Big(1+
\int_{0}^{t}\|\nabla v(\tau)\|_{L^\infty}d\tau\Big),
$$
where $C$ depends only on the dimension $d$ but not on the viscosity $\nu.$
\end{prop}
\begin{proof}
We denote by $\tilde{a}_{q}$ the unique global solution of the initial value problem:
$$\left\lbrace
\begin{array}{l}
\partial_t \tilde {a}_{q}+v\cdot\nabla \tilde{a}_{q}- \Delta \tilde{a}_{q}=\Delta_{q}f:=f_{q}
\\
{\tilde{a}_{q}}{(0)}=\Delta_{q}a^{0}.\\
\end{array}
\right.
$$
Using Proposition  \ref{propagation} with $r=+\infty$ and $s=\pm\frac12,$  one obtains
$$
\|\tilde{a}_{q}(t)\|_{B_{p,\infty}^{\pm\frac12}}\lesssim
\Big(\|\Delta_{q}{a}^0\|_{B_{p,\infty}^{\pm\frac12}}+
\int_{0}^t\|f_{q}(\tau)\|_{B_{p,\infty}^{\pm\frac12}}d\tau\Big) e^{C
\int_{0}^{t}\|\nabla v(\tau)\|_{{L^\infty}}d\tau}.
$$
Thus we deduce  from the definition of Besov spaces that for all $j\geq-1$
\begin{equation}\label{t1}
\|\Delta_{j}\tilde {a}_{q}\|_{L^\infty_{t}L^p}\lesssim 2^{-\frac12|q-j|}\Big(\|\Delta_{q}{a}^0\|_{L^p}+
\|f_{q}\|_{L^1_{t}L^p}\Big) e^{V(t)},
\end{equation}
with $V(t):= {C
\int_{0}^{t}\|\nabla v(\tau)\|_{L^\infty}d\tau}.$ 
Now by linearity one can write
$$
a(t,x)=\sum_{q\geq-1}\tilde {a}_{q}(t,x).$$
Taking  $N\in\NN$ that will  be carefully chosen later.  Then  we write by definition
\begin{eqnarray}\label{t2}
\nonumber\|a\|_{\widetilde L^\infty_{t}B_{p,r}^0}&\leq&
\Big(\sum_{j}\Big(\sum_{q}\|\Delta_{j}\tilde{a}_{q}\|_{L^\infty_{t}L^p}\Big)^r\Big)^{\frac{1}{r}}
\\
\nonumber&\leq&\Big(\sum_{j}\Big(\sum_{|q-j|\geq N}\|\Delta_{j}\tilde{a}_{q}\|_{L^\infty_{t}L^p}\Big)^r\Big)^{\frac{1}{r}}
+\Big(\sum_{j}\Big(\sum_{|q-j|< N}\|\Delta_{j}\tilde{a}_{q}\|_{L^\infty_{t}L^p}\Big)^r\Big)^{\frac{1}{r}}\\
&=&\hbox{I}+\hbox{II}.
\end{eqnarray}
To estimate the first term we use (\ref{t1}) and  the convolution inequality 
\begin{eqnarray}\label{t3}
\nonumber\hbox{I}&\lesssim& 2^{-\frac12 N}e^{V(t)}\|\big(\|{a_{q}}^0\|_{L^p}+
\|f_{q}\|_{L^1_{t}L^p}\big)_{q}\|_{\ell^r}\\
&\lesssim&2^{-\frac12 N}e^{V(t)}\Big(\|a ^0\|_{B_{p,r}^0}+\|f\|_{\widetilde L^1_{t}B_{p,r}^0}\Big).
\end{eqnarray}
To treat the second term of the right-hand  side of (\ref{t2}), we use two facts: the first
 one is that the operator $\Delta_{j}$ maps uniformly  $L^p$ into itself while the second  is  the $L^p$ energy estimate. So we find
\begin{eqnarray}\label{t4}
\nonumber\\
\nonumber\hbox{II}&\lesssim &\Big(\sum_{j}\Big(\sum_{|q-j|< N}\|\tilde{a}_{q}\|_{L^\infty_{t}L^p}\Big)^r\Big)^{\frac{1}{r}}\\
\nonumber&\lesssim&\Big(\sum_{j}\Big(\sum_{|q-j|< N}\|{a}_{q}^0\|_{L^p}+\|f_{q}\|_{L^1_{t}L^p}\Big)^r\Big)^{\frac{1}{r}}
\\
&\lesssim&N\big(\|a ^0\|_{B_{p,r}^0}+\|f\|_{\widetilde L^1_{t}B_{p,r}^0}\big).
\end{eqnarray}
Plugging estimates (\ref{t3}) and (\ref{t4}) into (\ref{t2}), we have
$$
\|a\|_{\widetilde L^\infty_{t}B_{p,r}^0}\lesssim\big(\|a ^0\|_{B_{p,r}^0}+
\|f\|_{\widetilde L^1_{t}B_{p,r}^0}\big)\big(
2^{-\frac12 N}e^{V(t)}+N\big).
$$
Taking
$$
N=\Big[\frac{2
V(t)}{ \log 2}+1\Big],
$$
leads to the desired inequality.
\end{proof}
%%%%%%
%%%%%
%%%%
%%%
%%
%
%\appendix
\section*{Appendix B.   Commutator estimate}
Our task now  is to prove the following commutator result.
 \begin{prop}\label{lemm2}
Let $u$ be a smooth function and  $v$ be a divergence-free vector field of $\RR^d$ such that its vorticity $\omega:=\hbox{curl }v$ belongs to
$L^\infty.$ Then we have for all  $q\geq -1,$
$$\big\|[\Delta_q, v\cdot \nabla ]u\big\|_{L^\infty}\lesssim
\|u\|_{L^\infty}\Big(\|\nabla\Delta_{-1}v\|_{L^\infty}+(q+2)\|\omega\|_{L^\infty}  \Big).$$
\end{prop}
 \begin{proof}
The principal tool is 
Bony's \mbox{decomposition \cite{b}:}
\begin{equation}\label{USA}[\Delta_q,v\cdot\nabla]u=[\Delta_q,T_v\cdot\nabla]u+[\Delta_q,
T_{\nabla\cdot}\cdot v]u+[\Delta_q,R(v\cdot\nabla,.)]u,
\end{equation}
where
\begin{eqnarray*}
&&[\Delta_q,T_v\cdot\nabla]u =\Delta_{q}(T_{v}\cdot\nabla u)-T_{v}\cdot\nabla\Delta_{q}u\\
&&[\Delta_q,T_{\nabla\cdot}\cdot v]u=\Delta_{q}(T_{\nabla u}\cdot v)-T_{\nabla\Delta_{q}u}\cdot v\\
&&
[\Delta_q,R(v\cdot\nabla,.)]u=\Delta_{q}(R(v\cdot\nabla, u))-R(v\cdot\nabla,\Delta_{q}u).
\end{eqnarray*}
From the definition of the paraproduct and according to Bernstein inequalities
\begin{eqnarray}\label{e7}
\nonumber\|[\Delta_q,
T_{\nabla\cdot}\cdot v]u\|_{L^\infty}&\lesssim& \sum_{|j-q|\leq
  4}\|S_{j-1}\nabla u\|_{L^{\infty}}\|\Delta_{j}v\|_{L^\infty}\\
  &\lesssim&\|u\|_{L^\infty}\|\omega\|_{L^\infty},
\end{eqnarray}
where  we have used here the following equivalence: $\forall j\in\NN,$
$$
\|\Delta_{j}v\|_{L^\infty}\approx 2^{-j}\|\Delta_{j}\omega\|_{L^\infty}.
$$
For the second term of  the right-hand side of (\ref{USA}), we have \begin{eqnarray*}
[\Delta_q,T_v\cdot\nabla]u&=&\sum_{j\geq 1}[\Delta_{q},S_{j-1}v\cdot\nabla\Delta_{j}]u,\\
&=&\sum_{|j-q|\leq 4}[\Delta_{q},S_{j-1}v\cdot\nabla]\Delta_{j}u.\\
\end{eqnarray*}
 To estimate each commutator, we write $\Delta_{q}$ as a convolution
\begin{eqnarray}
\nonumber
[\Delta_{q},S_{j-1}v\cdot\nabla]\Delta_{j}u(x)=2^{qd}\int
h(2^q(x-y))\big(S_{j-1}v(y)-S_{j-1}v(x)\big)\cdot\nabla\Delta_{j}
u(y)dy.
\end{eqnarray}
Thus, Young and Bernstein inequalities yield, for $|j-q|\leq 4,$ 
\begin{eqnarray}
\label{Tr1}
\qquad\big\|[\Delta_{q},S_{j-1}v\cdot\nabla]\Delta_{j}u\big\|_{L^\infty}&
\lesssim& 2^{-q}\|\nabla S_{j-1}v\|_{L^{\infty}}\|\Delta_{j}\nabla u
\|_{L^\infty}\\
\nonumber&\lesssim& \|\nabla S_{j-1}v\|_{L^\infty}\|u\|_{L^\infty}\\
\nonumber&\lesssim&\Big(\|\nabla\Delta_{-1}v\|_{L^\infty}+(q+2)\|\omega\|_{L^\infty} \Big)\|u\|_{L^\infty}.
\end{eqnarray}
Let us move to the remainder term. It can be written, in view of the
 definition, as
$$
\hbox{J}_{q}:=[\Delta_q,{R}(v\cdot\nabla,.)]u=\sum_{j\geq q-4,j\geq 0\atop\\
i\in\{\mp1,0\}}[\Delta_q,\Delta_{j}v]\cdot\nabla\Delta_{j+i} u
+\sum_{i\in\{0,1\}}[\Delta_q,\Delta_{-1}v]\cdot\nabla\Delta_{-1+i} u.
$$
It follows from the zero divergence condition that
\begin{equation*}
\hbox{J}_{q}=\sum_{i\in\{0,1\}}
[\Delta_q,\Delta_{-1}v]\cdot\nabla\Delta_{-1+i} u+\sum_{j\geq q-4,
 j\geq 0\atop\\i\in\{\mp1,0\}}\textnormal{div }\big( [
\Delta_q,\Delta_{j}v]\otimes\Delta_{j+i} u\big)\\
=\hbox{J}_{q}^1+\hbox{J}_{q}^2.
\end{equation*}
By the same way as (\ref{Tr1}) one has
\begin{eqnarray*}
\|\hbox{J}_{q}^1\|_{L^2}&\lesssim&2^{-q}\|\nabla\Delta_{-1}v\|_{L^\infty}\sum_{i=0}^1\|\nabla \Delta_{-1+i}u\|_{L^\infty}\\
&\lesssim&\|\nabla\Delta_{-1}v\|_{L^\infty}\| u\|_{L^\infty}.
\end{eqnarray*}
To estimate the second term we use Bernstein inequality
\begin{eqnarray*}\|\hbox{J}_{q}^2\|_{L^2}&\lesssim& \sum_{ j\geq q-4,j\geq 0
\atop
i\in\{\mp1,0\}}2^{q}\|\Delta_{j}
v\|_{L^{\infty}}\|\Delta_{j+i}u\|_{L^\infty}\\
&\lesssim &\|u\|_{L^\infty}\sum_{j\geq q-4}2^{q-j}\|\Delta_{j}\omega\|_{L^\infty}\\
&\lesssim&\|\omega\|_{L^\infty}\|u\|_{L^\infty},
\end{eqnarray*}

This completes the proof of  Proposition \ref{lemm2}.
\end{proof}
 %  \section{Appendix}

\end{document}